\documentclass[12pt,reqno]{amsart}
\usepackage{amsmath}
\usepackage{amssymb}

\newtheorem{thm}{Theorem}[section]
\newtheorem{lem}{Lemma}[section]

\newtheorem{coro}{Corollary}[section]\numberwithin{equation}{section}

\newfont{\bb}{msbm10 at 12pt}

\def\lf{\left}
\def\ri{\right}
\def\la{\langle}
\def\ra{\rangle}

\def\R{\hbox{\bb R}}

\def\V{\hbox{\bb V}}

\def\e{\epsilon}

\def\p{\partial}

\def\V{\mathcal{V}}
\def\mwy{\mathfrak{m}_{_{WY}}}
\def\mly{\mathfrak{m}_{_{ \rm LY}}}
\def\mby{\mathfrak{m}_{_{\rm BY}}}
\def\W{\mathcal{W}}
\def\madm{\mathfrak{m}_{_{ADM}}}

\def\tr{{\rm tr}}

\newcommand{\bal}{\begin{aligned}}      \newcommand{\eal}{\end{aligned}}
\newcommand{\ba}{\begin{array}}      \newcommand{\ea}{\end{array}}
\newcommand{\bc}{\begin{center}}     \newcommand{\ec}{\end{center}}
\newcommand{\be}{\begin{equation}}  \newcommand{\ee}{\end{equation}}
\newcommand{\beq}{\begin{eqnarray}}  \newcommand{\eeq}{\end{eqnarray}}
\newcommand{\beQ}{\begin{eqnarray*}} \newcommand{\eeQ}{\end{eqnarray*}}
\newcommand{\bi}{\begin{itemize}}    \newcommand{\ei}{\end{itemize}}
\newcommand{\bt}{\begin{tabular}}    \newcommand{\et}{\end{tabular}}
\newcommand{\bdm}{\begin{displaymath}} \newcommand{\edm}{\end{displaymath}}

\begin{document}
\title[]{Some estimates of Wang-Yau quasilocal energy}

\author{Pengzi Miao$^1$, Luen-Fai Tam$^2$ and Naqing Xie$^3$}

\address{School of Mathematical Sciences, Monash University,
Victoria, 3800, Australia.}
 \email{Pengzi.Miao@sci.monash.edu.au}

\address{The Institute of Mathematical Sciences and Department of
 Mathematics, The Chinese University of Hong Kong,
Shatin, Hong Kong, China.} \email{lftam@math.cuhk.edu.hk}

\address{School of Mathematical Sciences, Fudan
University, Shanghai 200433, China}
\email{nqxie@fudan.edu.cn}

\thanks{$^1$ Research partially supported by Australian Research Council Discovery Grant  \#DP0987650}

\thanks{$^2$Research partially supported by Hong Kong RGC General Research Fund  \#GRF 2160357}

\thanks {$^3$Research partially supported by the National Natural Science
Foundation of China \#10801036}

\renewcommand{\subjclassname}{%
  \textup{2000} Mathematics Subject Classification}
\date{}

\begin{abstract}
Given a spacelike $2$-surface $ \Sigma $ in a spacetime $ N$
and a constant future timelike unit vector $T_0 $ in $ \R^{3,1}$,
we derive  upper and lower estimates of  Wang-Yau quasilocal energy
$ E(\Sigma, X, T_0) $ for a given isometric embedding
$ X $ of $ \Sigma $ into a flat $ 3$-slice in $ \R^{3,1}$.
The quantity $ E(\Sigma, X, T_0) $ itself depends on the choice
of $X$, however the infimum of
 $ E( \Sigma, X ,T_0)$ over $ T_0 $ does not.
In particular,  when $ \Sigma $ lies in a time symmetric $3$-slice in $N$
and has nonnegative Brown-York quasilocal mass $ \mby(\Sigma)$,
our estimates show that $ \inf\limits_{T_0}E( \Sigma, X ,T_0) $
equals $ \mby (\Sigma)$.
We also study the spatial limit of $ \inf\limits_{T_0}E(S_r,X_r,T_0)$, where
 $ S_r $ is a large coordinate sphere in a fixed end of
an asymptotically flat initial data set $(M, g, p)$ and $ X_r $
is an isometric embeddings of  $S_r$ into
$\mathbb{R}^3 \subset \mathbb{R}^{3,1}$.
We show that if $(M, g, p)$ has future timelike ADM energy-momentum,
then $\lim\limits_{r\rightarrow\infty}\inf\limits_{T_0}E(S_r,X_r,T_0)$
equals  the ADM mass of $(M, g, p)$.
\end{abstract}

\maketitle \markboth{Pengzi Miao, Luen-Fai Tam and Naqing Xie}
{Some estimates of Wang-Yau quasilocal energy}


\section{Introduction} \label{intro}
In \cite{wy1,wy2}, Wang and Yau define a {\it quasilocal energy} $E(\Sigma, X,T_0)$ for a spacelike $2$-surface $\Sigma$ in a spacetime $N$, where $X$ is an isometric embedding of $\Sigma$ into the Minkowski spacetime $\R^{3,1}$ and $T_0$ is a constant future timelike  unit vector in $ \R^{3,1}$. Under the assumptions that the mean curvature vector $H$ of $\Sigma$ in $N$ and the mean curvature vector $H_0$ of
$ \Sigma $ when embedded in $\R^{3,1}$ are both spacelike, $E(\Sigma, X,T_0)$ is defined and can be expressed as follows (see \cite[(1.2),(1.5)]{wy3}):

\begin{equation}\label{qle}
\begin{split}E &(\Sigma, X, T_0)\\
&=\frac{1}{8\pi} \int_\Sigma \bigg\{ \sqrt{|{H}_0|^2(1+|\nabla\tau|^2)+(\Delta\tau)^2}-\sqrt{|{H}|^2(1+|\nabla\tau|^2)+(\Delta\tau)^2}\\
&\ \ \ -\Delta \tau \left[\sinh^{-1}(\frac{\Delta\tau}{\sqrt{1+|\nabla\tau|^2}|{H}_0|})-\sinh^{-1}(\frac{\Delta\tau}{\sqrt{1+|\nabla\tau|^2}|{H}|})\right]\\
&\ \ \ -\langle\nabla^{\R^{3,1}}_{\nabla\tau}
\frac{{J}_0}{|{H}_0|}, \frac{{H}_0}{|{H}_0|}\rangle
+\langle\nabla^{N}_{\nabla\tau} \frac{{J}}{|{H}|},
\frac{{H}}{|{H}|}\rangle\bigg\}
dv_\Sigma.\end{split}\end{equation}
Here $|H_0|=\sqrt{\la H_0,  H_0\ra}$, $|H|=\sqrt{\la H,H\ra }$, $\tau=-\la X,T_0\ra$,
$\nabla \tau$ and $\Delta \tau$ are the intrinsic gradient and the intrinsic Laplacian
of $ \tau $ on $\Sigma$;
$J$ is the future timelike unit normal vector
field along $\Sigma$ in $N$ which is dual to $H$ along the
light cone in the normal bundle of $ \Sigma $. Namely, if $e_1$, $e_2$,
$e_3$, $e_4$  are orthonormal tangent vectors of $N $ such
that $e_1, e_2$ are tangent to $\Sigma$, $e_3$ is spacelike
with $ \la e_3, H \ra < 0 $ and $e_4$ is future timelike,
then $H=\la H,e_3\ra e_3- \la H,e_4\ra e_4$ and $J=\la
H,e_4\ra e_3- \la H,e_3\ra e_4$; $J_0$ is defined similarly
for $\Sigma$ when embedded in $\R^{3,1}$.

Wang and Yau \cite{wy1,wy2} define a new {\it quasilocal
mass} of $ \Sigma $ in $ N$, which we denote by
$\mwy (\Sigma)$,  to be the
infimum of $E(\Sigma,X,T_0)$ over those $ X $ and $ T_0 $
such that the resulting $ \tau = - \la T_0, X \ra $ is
{\it admissible} (see \cite[Definition 5.1]{wy2} for
the definition of admissible data). They prove an important
property on positivity of $\mwy (\Sigma)$ under the
assumption that $N$ satisfies the usual dominant energy condition.
Note that, if $ \Sigma $ has positive Gaussian curvature, then
$ \tau = 0 $ is admissible \cite[Remark 1.1]{wy2}. Hence,
$ \mwy( \Sigma ) \leq E ( \Sigma, {X}, \hat{T}_0 ) =
\mly ( \Sigma )$,
where $ {X} $ is an embedding of $ \Sigma $ into
$\R^3=\{(0,x)\in \R^{3,1}\}$, $ \hat{T}_0 = ( 1, 0, 0, 0)$,
and $ \mly (\Sigma) $ is the Liu-Yau quasilocal mass of
$ \Sigma $ \cite{ly1, ly2}.
Note that $ \mly (\Sigma) $ equals the Brown-York quasilocal mass
$\mby(\Sigma)$  \cite{by1, by2} if in addition $\Sigma$ bounds
a compact time-symmetric hypersurface.

The expression of $ E( \Sigma, X, T_0)$
is rather complicated.
It is not clear if $ \mwy( \Sigma) $ can be achieved by
some admissible data.
In \cite{wy3}, restricting to an embedding $ X $ of $ \Sigma$
into $ \R^3 \subset \R^{3,1}$,
 Wang and Yau study the spatial limit of $ E( \Sigma, X, T_0) $
on an asymptotically flat spacelike hypersurface.
Motivated by their work,   we want to get some lower and upper
estimates of $ E( \Sigma, X, T_0) $ for a fixed $ 2$-surface $ \Sigma$.


More precisely, let us assume that $ \Sigma $ has positive Gaussian curvature.
Isometrically embed $ \Sigma $ in  $\R^3=\{(0,x)\in \R^{3,1}\}$
and let $ X $ be the  embedding.
 Let $T_0=(\sqrt{1+|a|^2},a)$, where $a=(a^1,a^2,a^3)\in \R^3$, be a constant future timelike unit vector in  $\R^{3,1}$.   In this work, we will prove the following estimates:
{\it
\begin{equation}\label{mainresult-e1}
\begin{split}
\sqrt{ 1 + | a |^2
}\lf(\mly (\Sigma)+C\ri) - &\sum_{i=1}^3 a^i \V_i
 \\
 \ge& E ( \Sigma, X, T_0 )\\
\geq& \sqrt{ 1 + | a |^2
} \mly (\Sigma) - \sum_{i=1}^3 a^i \V_i\\
 \end{split}
\end{equation}
where $C$ is a constant depending only on $|H_0|$, $| H |$
and their integrals on $ \Sigma$,
and $\V=(\V_1,\V_2,\V_3)$ is a constant vector in $ \R^3$
so that $ - \sum_{i=1}^3 a^i \V_i$ equals
$  \int_\Sigma \langle\nabla^{N}_{\nabla\tau} \frac{{J}}{|{H}|},
\frac{{H}}{|{H}|}\rangle d v_\Sigma$
in \eqref{qle}}.
Note that, if we let
$$
\W=(\mly (\Sigma),\V),
$$
then $\sqrt{ 1 + | a |^2 } \mly (\Sigma) - \sum_{i=1}^3 a^i \V_i$ is exactly $-\la T_0,\W\ra$, and \eqref{mainresult-e1} can be written as:
\begin{equation}\label{mainresult-e1-1}
\begin{split}
-\la T_0,\W\ra+C\sqrt{ 1 + | a |^2
}
 \ge  E ( \Sigma, X, T_0 )\ge
-\la T_0,\W\ra.
 \end{split}
\end{equation}

An immediate application of the estimates \eqref{mainresult-e1}
is the following result relating the Wang-Yau quasilocal energy
 and the Brown-York quasilocal mass:

{\it With the above notations, suppose $ \Sigma $ bounds a compact, time-symmetric
hypersurface $ \Omega $ in a spacetime $ N $ satisfying the
dominant energy condition. Suppose $ \Sigma $ has positive
Gaussian curvature and has positive
mean curvature  in $ \Omega $ with respect to the
outward unit normal, then
$$ \inf_{T_0} E( \Sigma, X, T_0 )  =  \mby ( \Sigma ) . $$}

Another application of  \eqref{mainresult-e1}  is in the study
of  the spatial limit of Wang-Yau quasilocal energy
on an asymptotically flat spacelike hypersurface.
Recall that a spacelike hypersurface $M$ with
induced metric $g$ and second fundamental form $p$ in a
spacetime $N$ is called {\it asymptotically flat}   if there is
a compact set $K$ such that $M\setminus K$  has finitely
many ends, each of which is diffeomorphic to the complement
of a Euclidean ball in  $\R^3$ and such that the metric $g$
is of the form $g_{ij}=\delta_{ij}+a_{ij}$ so that
 \begin{equation}\label{af-e1}
 r|a_{ij}|+r^{2}|\p a_{ij}|+r^{3}|\p\p a_{ij}| \le C
 \end{equation}
 and $p_{ij}$ satisfies
\begin{equation}\label{af-e2}
 r^2|p_{ij}|+r^{3}|\p p_{ij}| \le C
\end{equation}
for some constant $C$, where $r=|x|$ denotes the coordinate
length in $ \R^3$.  In \cite{wy3}, Wang and Yau prove the following:
\begin{thm}\cite[Theorem 3.1]{wy3} \label{WY}
Let $(M,g,p)$ be an asymptotically flat spacelike hypersurface in a spacetime $N$.
On any given end of $ M$, let $S_r$ be the coordinate sphere of radius $ r $,
let $X_r$ be the isometric embedding of $S_r$
into $\R^3=\{(0,x)\in \R^{3,1}\}$ given by Lemma \ref{fst-l}, and let
$T_0=(\sqrt{1+|a|^2},a^1,a^2,a^3)$ be a fixed future timelike unit
vector in $\R^{3,1}$. Then
\beq \label{wylimit}
\lim_{r\rightarrow\infty}E(S_r,X_r,T_0)=\sqrt{1+|a|^2}E+\sum_{i=1}^3 a^iP_i ,
\eeq
where $E$ is the ADM energy of $ (M, g, p)$ and $P = ( P_1, P_2, P_3)$
is the ADM linear momentum of $(M, g, p)$.
\end{thm}
Since the Wang-Yau quasilocal mass $ \mwy(\Sigma) $
is defined as the infimum of a  class of $E(\Sigma, X, T_0)$,
we would like to understand the asymptotical behavior of the
infimum   $\inf\limits_{T_0}E(S_r,X_r,T_0)$   as $r\to \infty$.
Note that in Theorem \ref{WY}, if $E>|P|$ and hence the ADM mass \cite{adm} $\madm $
of $(M, g, p)$ is positive, then it is not hard to see that
 \be \label{infadm}
 \inf_{a\in \R^3}\lf(\sqrt{1+|a|^2}E+\sum_{i=1}^3 a^iP_i\ri)= \sqrt{E^2-|P|^2}= \madm .
 \ee
Applying the estimates \eqref{mainresult-e1}, we will prove the following:

{\it With the notations given in Theorem \ref{WY}, suppose $N$ satisfies
the dominant energy condition and suppose $N$ is not flat along $M$, then
\begin{equation}\label{mainresult-e2}
    \lim_{r\to\infty}\inf_{T_0}E(S_r,X_r,T_0)=\madm  .
\end{equation}}

Again, since $ \mwy (S_r) $ is defined to be the infimum of $ E(S_r, X_r, T_0) $ over those isometric embeddings $ X_r:\Sigma\hookrightarrow \R^{3,1} $ and
 constant future timelike unit vectors $ T_0 $ in $\R^{3,1}$
such that $ \tau = - \la X_r, T_0 \ra $ are admissible, we remark that
it is yet unclear
 whether $\lim\limits_{r\to\infty} \mwy (S_r)= \madm $.

\section{Estimates of Wang-Yau quasi-local energy}

In this section, we derive the main estimates of
$E(\Sigma,X, T_0)$ for a given
isometric embedding $ X $ of $ \Sigma $ into a  flat $ 3$-slice in $ \R^{3,1}$.
Precisely, let $ \Sigma $ be a  spacelike $2$-surface
in a spacetime $ N$ such that $ \Sigma $ has positive
Gaussian curvature and the mean curvature vector $ H $ of $
\Sigma $ in $ N $ is spacelike. Let  $ X $ be an isometric
embedding of $ \Sigma $ into $ \R^3 = \{ (0, x) \in
\R^{3,1}  \}$ and let $ T_0 \in \R^{3,1} $ be an arbitrary constant
future timelike unit vector. Let $ E( \Sigma, X, T_0 )$ be
given by \eqref{qle} with  $ \tau = - \la T_0, X  \ra $.
Since $X$ is only unique up to a rigid motion of $\R^3$,
we note that $E(\Sigma, X,T_0)$  depends on the choice of $X$,
but the infimum
\be \label{eq-inf}
  \inf_{ T_0 } E ( \Sigma, X, T_0 )
\ee
does not.  Hence, if it is finite,  \eqref{eq-inf} gives a geometric invariant of $ \Sigma $ in $ N$.

First, we define a constant vector $ \W = \W(\Sigma,X) \in \R^{3,1}$
associated to  $ \Sigma $ and $ X$. Let $ \alpha_{ \hat{e}_3
} ( \cdot ) $ be the connection $ 1$-form on $ \Sigma $
introduced in \cite[(1.3)]{wy2} which is  defined as
\begin{equation} \label{c1form}
\alpha_{\hat{e}_3} ( Y ) =  \la \nabla^N_Y \frac{J}{|H|} ,
\frac{ H }{ | H |} \ra , \ \ \forall \ Y \in T \Sigma  .
\end{equation}
Let $ V=V(\Sigma) $ be the vector field on $ \Sigma $
that is dual to $ \alpha_{\hat{e}_3} ( \cdot)  $. Given an
embedding $ X : \Sigma \hookrightarrow \R^3 $, we identify $ V
$ with $ d X ( V ) $ through the tangent map $ d X $ and
hence view $ V = V( \Sigma, X) $ as an $ \R^3$-valued vector field along $
\Sigma $.  Now define
\begin{equation} \label{dfofV}
\V (\Sigma,X)= \frac{1}{8\pi} \int_\Sigma V d v_{\Sigma}  \in \R^3
\end{equation}
and
\begin{equation}\label{dfofW}
\W(\Sigma,X) = ( \mly ( \Sigma) , \V(\Sigma,X) )  \in \R^{3,1} .
\end{equation}
Here $ \mly  (\Sigma )$ is  the Liu-Yau quasilocal mass
defined in \cite{ly1, ly2}, i.e.
\begin{equation}
\mly ( \Sigma ) = \frac{ 1 }{ 8 \pi} \int_\Sigma ( k_0 - |
H | ) d v_\Sigma,
\end{equation}
where $ k_0 $ is the mean curvature of $ \Sigma $ with
respect to the outward unit normal when isometrically
embedded in $ \R^3$. Note that if $ \Sigma $ bounds a
compact, time-symmetric hypersurface $ \Omega $ in $ N $
and $ \Sigma $ has positive mean curvature  in $ \Omega$
with respect to the outward unit normal, then $ \mly
(\Sigma )$ agrees with the Brown-York quasilocal mass $
\mby (\Sigma )$  \cite{by1, by2}.

If $ X $ differs by a rigid motion in $ \R^{3}$, then $\V(\Sigma,X)
$ differs by a corresponding rotation in $ \R^3$ and $\W (\Sigma,X)$
differs by a rotation in $\R^3$ considered as a Lorentzian
transformation in $\R^{3,1}$.

\begin{thm} \label{estinf}
Let $ \Sigma$, $ X $, $ T_0 $ and $ \V=\V(\Sigma,X) $ and
$\W=\W(\Sigma,X)$ be given as above. The followings are true:

\begin{enumerate}
\item[(i)] Suppose  $ \V = ( \V_1, \V_2, \V_3  )$, then
\begin{equation} \label{mainest1}
\begin{split}
-\la T_0,\W\ra+C\sqrt{ 1 + | a |^2 }=&  \sqrt{ 1 + | a |^2 }
\lf(\mly (\Sigma)+C\ri) - \sum_{i=1}^3 a^i \V_i \\
\ge& E ( \Sigma, X, T_0 ) \\
\geq &\sqrt{ 1 + | a |^2
} \mly (\Sigma) - \sum_{i=1}^3 a^i \V_i\\
 =&-\la T_0,\W\ra,
 \end{split}
\end{equation}
where $ T_0=(\sqrt{1+|a|^2},a^1,a^2,a^3)$ and
 $C$ is the constant given by
$$
C=\sup_{\Sigma}
\lf|\lf(\frac{|H_0|^2}{|H|^2}+\frac{|H_0|}{|H|}-2\ri)\ri|
\frac{1}{8\pi} \int_\Sigma\lf|\lf(|H_0|-|H|\ri)\ri| .
$$
Moreover, the equality
\be \label{eqcase}
 E ( \Sigma, X, T_0 ) = -\la T_0,\W\ra
\ee
holds if and only if $ T_0 = (1, 0, 0, 0) $ or
$ | H | = | H_0 | $
everywhere on $ \Sigma $.
In this second case, $ E(\Sigma, X, T_0 ) = - \sum_{i=1}^3 a^i \V_i $,
$ \forall$ $ T_0$.

\item[(ii)] If $\W$ is future timelike, then
\begin{equation}\label{inf-e1}
\begin{split}
\sqrt{-\la\W,\W\ra} + C \frac{ \mly(\Sigma) }{ \sqrt{-\la\W,\W\ra}  } \ge  & E(\Sigma,X,T_0^*)  \\
\ge &
\inf_{T_0}E(\Sigma,X,T_0)
\ge  \sqrt{-\la\W,\W\ra} ,
\end{split}
\end{equation}
where $ T_0^* = \frac{\W}{\sqrt{-\la\W,\W\ra}}$.

\item[(iii)]  The vector $ - \V $ is the gradient of the function $ E ( a)
= E ( \Sigma, X, T_0) $ at $ a = (0, 0, 0)$.
If in addition $ \mly ( \Sigma ) \geq 0 $, then $ \V = 0 $ if and only if
\be \label{infmly}
 \inf_{ T_0} E ( \Sigma, X, T_0 ) = \mly ( \Sigma).
 \ee

\end{enumerate}

\end{thm}

\begin{proof}
(i) Let $ \{ e_\alpha \ | \ \alpha = 1, 2 \}$ be an orthonormal
frame of $ T_p \Sigma $ for any $ p \in \Sigma$.
Identifying $ e_\alpha $ with $ d X ( e_\alpha )$ through
the tangent map $ d X$, we have
\begin{equation*}
\nabla \tau =  - \sum_{\alpha = 1 }^2  \lf(  e_\alpha \la
a, X \ra \ri) e_\alpha =  - \sum_{\alpha = 1 }^2  \la a,
e_\alpha  \ra  e_\alpha .
\end{equation*}
Hence,
\begin{equation} \label{cinN}
\langle\nabla^{N}_{\nabla\tau} \frac{{J}}{|{H}|},
\frac{{H}}{|{H}|}\rangle =   - \sum_{\alpha = 1 }^2  \la a,
e_\alpha  \ra  \langle\nabla^{N}_{ e_\alpha }
\frac{{J}}{|{H}|}, \frac{{H}}{|{H}|}\rangle = - \la a, V
\ra,
\end{equation}
where $ V $ is the vector field on $ \Sigma $ dual to the
connection $ 1$-form $ \alpha_{\hat{e}_3} (\cdot) $ defined
in \eqref{c1form}. On the other hand,
 \begin{equation} \label{cinR}
 -\langle\nabla^{\R^{3,1}}_{\nabla\tau} \frac{{J}_0}{|{H}_0|}, \frac{{H}_0}{|{H}_0|}\rangle  = 0
 \end{equation}
since $X ( \Sigma )$ lies in $\R^3$. Therefore, by
\eqref{qle}, \eqref{cinN}-\eqref{cinR}, we have
\begin{equation} \label{qle1}
E (\Sigma, X, T_0) = \tilde{E} (\Sigma, X, T_0)  - \la a,
\V \ra ,
\end{equation}
where
\begin{equation} \label{qletE}
\begin{split}
 & \tilde{E} (\Sigma, X, T_0)\\
= & \ \frac{1}{8\pi} \int_\Sigma \bigg\{ \sqrt{|{H}_0|^2(1+|\nabla\tau|^2)+(\Delta\tau)^2}-\sqrt{|{H}|^2(1+|\nabla\tau|^2)+(\Delta\tau)^2}\\
&  \ -\Delta \tau
\left[\sinh^{-1}(\frac{\Delta\tau}{\sqrt{1+|\nabla\tau|^2}|{H}_0|})-\sinh^{-1}(\frac{\Delta\tau}{\sqrt{1+|\nabla\tau|^2}|{H}|})\right]
\bigg\} d v_\Sigma .
\end{split}
\end{equation}
We first prove that:
\begin{equation} \label{mainest2}
 \tilde{E} (\Sigma, X, T_0)  \geq \sqrt{ 1 + | a |^2} \mly( \Sigma) .
\end{equation}

To prove \eqref{mainest2}, we use spherical coordinates on
$ \R^{3}$. Let $ \rho \geq 0 $ be a scalar  and  $ \omega
\in \mathbb{S}^2 $ be a unit vector in $ \R^3 $. Write $ a
= \rho \omega $, we have
$$
\tau = - \la a , X \ra = - \rho \la \omega, X \ra ,
$$
$$
\Delta \tau = - \la a, \Delta X \ra =  - \la \rho \omega,
H_0 \ra = \rho k_0 \la \omega, e^{ H_0} \ra ,
$$
$$
1 + | \nabla \tau |^2 =  1 + \rho^2 ( 1 -  \la \omega ,
e^{H_0} \ra^2 ) ,
$$
where $  H_0 $ is the mean curvature vector of $ X( \Sigma)
$ in $ \R^3 \subset \R^{3,1} $, $ e^{ H_0} = - \frac{ H_0}{
|  H_0 |}$ is the outward unit normal to $ X( \Sigma) $ in
$ \R^3$, and $ k_0 = | H_0 |  $ is the mean curvature of $
X( \Sigma) $ with respect to $ e^{H_0} $  in $ \R^3 $.

In what follows, we let
\be \label{dfofp}
 p = \la \omega, e^{H_0 } \ra , \ \
 q = \sqrt{ 1 -  \la \omega, e^{ H_0 } \ra^2 }
 \ee
be two functions on $ \Sigma $. In terms of $ p $, $ q$, we
have \be \Delta \tau = \rho k_0 p , \ \ 1 + | \nabla \tau
|^2 = 1 + \rho^2 q^2 \ee and \be \frac{ \Delta \tau}{
\sqrt{ 1 + | \nabla \tau |^2 } } = k_0 \frac{ \rho p}{
\sqrt{ 1 + \rho^2 q^2} } . \ee The quantity $ \tilde{E}
(\Sigma, X, T_0 )$ now becomes
 a function of  $ \rho $ and $ \omega$, say
 $ \tilde{E}( \rho, \omega)$.
 By \eqref{qletE}, we have
\begin{equation*}
\begin{split}
\tilde{E}( \rho, \omega)  = & \tilde{E} (\Sigma, X, T_0 ) \\
= & \frac{1}{ 8 \pi}  \int_\Sigma  \bigg\{ \lf[  \sqrt{
k_0^2 ( 1 + \rho^2 ) } -
\sqrt{ k^2_0 \rho^2 p^2 + k^2 ( 1 + \rho^2 q^2 ) }  \ri]  \\
&  + \rho k_0 p  \lf[ \sinh^{-1} \lf( \frac{ \rho p } {
\sqrt{ 1 + \rho^2 q^2 } }  \frac{ k_0} {k } \ri) -
\sinh^{-1} \lf( \frac{   \rho p} { \sqrt{ 1 + \rho^2 q^2 }
} \ri)  \ri] \bigg\}  ,
\end{split}
\end{equation*}
where $ k = |  H | > 0 $ and $ H $ is the mean curvature
vector of $ \Sigma $ in $ N $. Henceforth, we omit writing
the volume form $ d v_\Sigma $ for simplicity.

Since $ k_0 > 0 $, we can rewrite $ \tilde{E} ( \rho,
\omega)  $ as
\begin{equation} \label{defoftE}
\begin{split}
 \tilde{E} ( \rho, \omega)
 = & \frac{1}{ 8 \pi}  \int_\Sigma  k_0 \bigg\{  \lf[  \sqrt{ 1 + \rho^2  } -
\sqrt{  \rho^2 p^2 + t^2 ( 1 + \rho^2 q^2 ) }  \ri]  \\
& + k_0 (\rho p ) \lf[ \sinh^{-1} \lf( \frac{ \rho p } {
\sqrt{ 1 + \rho^2 q^2 } }  t^{-1} \ri) - \sinh^{-1} \lf(
\frac{   \rho p} { \sqrt{ 1 + \rho^2 q^2 } } \ri)  \ri]
\bigg\} ,
\end{split}
\end{equation}
where $ t = \frac{ k }{k_0} > 0  $ is a function on $
\Sigma$.  Let \be \label{doff}
 f ( \rho, \omega ) =  \frac{ \rho p} { \sqrt{ 1 + \rho^2 q^2} },
 \ee
then
 \be
  1+\rho^2q^2=\frac{1+\rho^2}{1+f^2}, \ \ p=\frac{f}{\rho}
  \frac{(1+\rho^2)^\frac12}{(1+f^2)^\frac12}.
  \ee
 Now define
  \be\begin{split}
         B( \rho, \omega ) = & \sqrt{ 1 + \rho^2  } - \sqrt{
\rho^2 p^2 + t^2 ( 1 + \rho^2 q^2 ) }  , \\
         =&
         \sqrt{1+\rho^2}\lf(1-\frac{\lf(t^2+f^2\ri)^\frac12}{\lf(1+f^2\ri)^\frac12}\ri)
     \end{split}
     \ee
and \be
 F ( \rho , \omega ) = \rho p \lf[ \sinh^{-1} \lf( \frac{f }{t} \ri) - \sinh^{-1} ( f ) \ri]  .
 \ee
By \eqref{defoftE}, we have
  \begin{equation}\label{E(rho)-e2}
 \tilde{E} ( \rho, \omega)= \frac{1}{ 8 \pi} \int_\Sigma k_0
 [ B( \rho, \omega)+F( \rho, \omega) ] .
  \end{equation}
Direct calculation gives \be  \label{pf} \frac{ \p f}{ \p
\rho } = p ( 1 +  \rho^2 q^2 )^{- \frac{3}{2}
}=\frac{f(1+f^2) }{\rho(1+\rho^2) },
 \ee
 \be
\begin{split} \label{pB}
\frac{ \p B }{ \p \rho}    = &   \frac{ \rho}{ \sqrt{ 1 +
\rho^2 }
}\lf(1-\frac{\lf(t^2+f^2\ri)^\frac12}{\lf(1+f^2\ri)^\frac12}\ri)\\
&-(1+\rho^2)^\frac12\lf(\frac{f\frac{\p f}{\p
\rho}}{\lf(t^2+f^2\ri)^\frac12\lf(1+f^2\ri)^\frac12}-\frac{f\frac{\p
f}{\p \rho}\lf(t^2+f^2\ri)^\frac12}{\lf(1+f^2\ri)^\frac32}
\ri)
 \\
 =& \frac{ \rho}{ ( 1 +\rho^2 )^\frac12 }-\frac{(1-t^2)f^2+\rho^2(t^2+f^2)}{\rho\lf(t^2+f^2\ri)^\frac12\lf(1+f^2\ri)^\frac12\lf(1+\rho^2\ri)^\frac12} ,
\end{split}
\ee \be
\begin{split}\label{pF}
\frac{\p F}{\p \rho } = &  p \lf[ \sinh^{-1} \lf( \frac{f }{t} \ri) - \sinh^{-1} ( f ) \ri] + \\
& \rho p \lf[ \lf( 1 + \frac{ f^2 }{ t^2 } \ri)^{-\frac 12}
\frac{1}{t} - \lf( 1 + f^2 \ri)^{-\frac 12}  \ri]
\frac{df}{d \rho }\\
=& -\frac{f}{\rho}
  \frac{(1+\rho^2)^\frac12}{(1+f^2)^\frac12}\lf[\sinh^{-1}f+\frac{f(1+f^2)^\frac12 }{(1+\rho^2) }\ri] \\
  &+\frac{f}{\rho}
  \frac{(1+\rho^2)^\frac12}{(1+f^2)^\frac12}\lf[\sinh^{-1}\lf(\frac ft\ri)+\frac{f(1+f^2) }{(1+\rho^2)(t^2+f^2)^\frac12 }\ri] ,
\end{split}
\ee where we assume $ \rho > 0 $ whenever $ \rho $ appears
in a denominator.
By \eqref{pB} and \eqref{pF}, we have
 \be \label{ptE-by}
 \begin{split}
& \frac{ \p \tilde{E} }{ \p \rho} -\frac{\rho}{(1+\rho^2)^\frac12} \mly (\Sigma)  \\
= & \ \frac{1}{8\pi}   \int_\Sigma k_0\lf[
\frac{\p B}{\p \rho}+\frac{\p F}{\p \rho}-\frac{\rho}{(1+\rho^2)^\frac12}(1-t)\ri] \\
= &\  \frac{1}{8\pi}   \int_\Sigma
k_0\lf(\Phi(1)-\Phi(t)\ri) ,
\end{split}
\ee where
 \be
 \begin{split}
 \Phi(t) =&\frac{(1-t^2)f^2+\rho^2(t^2+f^2)}{\rho\lf(t^2+f^2\ri)^\frac12\lf(1+f^2\ri)^\frac12\lf(1+\rho^2\ri)^\frac12} \\
 & -\frac{f}{\rho}
  \frac{(1+\rho^2)^\frac12}{(1+f^2)^\frac12}\lf[\sinh^{-1}\lf(\frac ft\ri)+\frac{f(1+f^2) }{(1+\rho^2)(t^2+f^2)^\frac12 }\ri]\\
  &-\frac{\rho}{(1+\rho^2)^\frac12}t  \\
  =&\frac{1}{\rho\lf(1+f^2\ri)^\frac12\lf(1+\rho^2\ri)^\frac12}
  \lf[-f(1+\rho^2)\sinh^{-1}\lf(\frac ft\ri)+
   (\rho^2-f^2) (t^2+f^2)^\frac12 \ri] \\
   &
   -\frac{\rho}{(1+\rho^2)^\frac12}t.
\end{split}
 \ee
Now
 \be\label{Phi-e1}
 \begin{split}
 \frac{\p \Phi}{\p t}=& \frac{f^2(1+\rho^2)+(\rho^2-f^2)t^2}{t\rho\lf(1+f^2\ri)^\frac12\lf(1+\rho^2\ri)^\frac12\lf(t^2+f^2\ri)^\frac12}
 -\frac{\rho}{(1+\rho^2)^\frac12}\\
  =& \lf[ \frac{(1-t^2)f^2 +\rho^2(t^2+f^2)}{t \rho^2 \lf(1+f^2\ri)^\frac12 \lf(t^2+f^2\ri)^\frac12}
 - 1 \ri]  \frac{ \rho }{(1+\rho^2)^\frac12} .
 \end{split}
 \ee
Hence if $ 0< t\le 1$,
\be \label{Phi-e1-2}
 \begin{split}
 \frac{ \p \Phi}{\p t}\ge& \frac{\rho}{(1+\rho^2)^\frac12}
 \lf[\frac{(t^2+f^2)^\frac12}{t(1+f^2)^\frac12}-1\ri]\\
 =&\frac{\rho}{(1+\rho^2)^\frac12}
 \lf[\frac{\lf(1+\lf(\frac ft\ri)^2\ri)^\frac12}{ (1+f^2)^\frac12}-1\ri]\\
 \ge&0.
 \end{split}
 \ee
 Similarly, if $t>1$, then
 \be \label{Phi-e1-3}
 \begin{split}
 \frac{ \p \Phi}{\p t}\le& \frac{\rho}{(1+\rho^2)^\frac12}
 \lf[\frac{(t^2+f^2)^\frac12}{t(1+f^2)^\frac12}-1\ri]\\
 =&\frac{\rho}{(1+\rho^2)^\frac12}
 \lf[\frac{\lf(1+\lf(\frac ft\ri)^2\ri)^\frac12}{ (1+f^2)^\frac12}-1\ri]\\
 \le&0.
 \end{split}
 \ee
Therefore, $\Phi(1)=\max_{t>0}\Phi(t)$.
By \eqref{ptE-by}, we have
\begin{equation} \label{maindest}
\frac{ \p \tilde{E}}{\p \rho} \ge
\frac{\rho}{(1+\rho^2)^\frac12} \mly (\Sigma)
\end{equation}
for all $ \rho > 0 $. Integrating \eqref{maindest} and note
that $ \tilde{E}(0, \omega) = \mly(\Sigma) $, we conclude
\be \label{mainest3}
\tilde{E}(\rho, \omega) \geq \sqrt{ 1 + \rho^2} \mly
(\Sigma) ,
\ee
which proves \eqref{mainest2}. Now the lower
bound of $E(\Sigma,X,T_0)$ in  \eqref{mainest1} follows
directly from \eqref{qle1} and \eqref{mainest2}.

In order to obtain an upper bound for $E(\Sigma,X,T_0)$, by
\eqref{Phi-e1} and using the fact $f^2\le \rho^2$, we have,
if $ 0 < t_0 \le t \le 1$,
\be
 \begin{split}
 \frac{ \p \Phi}{\p t} \le &
\lf( \frac{1}{t^2}+\frac{1}{t}-2\ri)
 \frac{\rho}{(1+\rho^2)^\frac12}  \\
 \le & \lf( \frac{1}{t_0^2}+\frac{1}{t_0}-2\ri)
 \frac{\rho}{(1+\rho^2)^\frac12}  .
 \end{split}
 \ee
Hence by the mean value theorem,  if $ 0 < t_0 \le 1 $,
\be\label{Phi-e2}
 \Phi(1)-\Phi(t_0)   \le   ( 1 - t_0)  \lf( \frac{1}{t_0^2}+\frac{1}{t_0} - 2 \ri)
 \frac{\rho}{(1+\rho^2)^\frac12 }  .
    \ee
Similarly, if $ t_0 \ge t \ge 1$, we have
\be
 \begin{split}
 \frac{ \p \Phi}{\p  t}  \ge& \lf(\frac{1}{t^2}+\frac{1}{t}-2\ri)
 \frac{\rho}{(1+\rho^2)^\frac12}  \\
 \ge & \lf(\frac{1}{t_0^2}+\frac{1}{t_0}-2\ri)
 \frac{\rho}{(1+\rho^2)^\frac12} .
 \end{split}
 \ee
Thus \eqref{Phi-e2} is also true for $ t_0 \ge 1$.
By \eqref{ptE-by} and \eqref{Phi-e2},
 \be \label{ptE-by-1}
 \begin{split}
& \frac{ \p \tilde{E} }{ \p \rho} -\frac{\rho}{(1+\rho^2)^\frac12} \mly (\Sigma)  \\
\le  & \frac{\rho}{(1+\rho^2)^\frac12}
 \sup_{\Sigma}
\lf|\lf(\frac{|H_0|^2}{|H|^2}+\frac{|H_0|}{|H|}-2\ri)\ri|
 \frac{1}{ 8 \pi } \int_\Sigma\lf|\lf(|H_0|-|H|\ri) \ri|  .
\end{split}\ee
From this, it is easy to see that the upper bound
for $E(\Sigma,X,T_0)$ in \eqref{mainest1} holds.

To prove the rest part of (i), we note that
if  $T_0 = (1, 0,0,0)$, then
$ E (\Sigma, X, T_0) = \mly(\Sigma) $;
 if  $|H|=|H_0|$, then $\tilde{E}(X,\Sigma,T_0)=0$.
Hence \eqref{eqcase} holds automatically in either case.
Now suppose \eqref{eqcase} is true for some
$T_0 $ with $  a  = \rho_0 \omega $ where
$ \rho_0 > 0 $, then by the   proof of  \eqref{mainest3},
we  have
$  \Phi (t) = \Phi (1) $
 everywhere on $ \Sigma $  for any $ 0 < \rho \le \rho_0$.
A detailed examination of \eqref{Phi-e1}-\eqref{Phi-e1-3}
 shows that, at points $ x \in \Sigma $ where
 $ f = 0 $, we have $ \frac{ \p \Phi}{\p t} = 0 $, $ \forall t $;
 while at points $ x \in \Sigma $ where $ f \neq 0$, we have
 $ \frac{ \p \Phi}{\p t} >  0 $, $ \forall  0 < t <1 $
 and
 $ \frac{ \p \Phi}{\p t} <  0 $, $ \forall  t > 1 $.
Therefore,
 the fact $ \Phi(t) = \Phi(1) $ on $ \Sigma$  implies that
$ t = 1 $ on the subset $ \{ f \neq 0 \}
\subset \Sigma $. On the other hand, since $ \rho_0 > 0 $,
we know $ f = 0 $ if and only if $ p = 0$ or equivalently
 $  \la \omega, e^{H_0 }  \ra = 0 $.  Since
 $ X(\Sigma) $ is a strictly convex close surface in $ \R^3$,
 the Gauss map that sends a point on $ \Sigma $ to
 its outward unit normal $ e^{H_0 }$ is
 a diffeomorphism from $ \Sigma$ to  $ \mathbb{S}^2$.
 Therefore, the set $ \{ \la \omega , e^{H_0 }  \ra = 0   \} $ is a closed embedded curve in $ \Sigma$. Consequently, its complement
 $ \{ f \neq 0 \}$  is dense in $ \Sigma$. Therefore,
 $ t = 1 $ and hence  $ | H | = | H_0 |$ everywhere on $ \Sigma$.
 In this case, by definition, $ \tilde{E}( \Sigma, X, T_0) = 0 $, and
 by \eqref{qle1}, $ E ( \Sigma, X, T_0 ) = - \sum_{i=1}^3 a^i \V_i $.

(ii) Suppose $\W$ is future timelike, then
$-\la T_0,\W \ra$ attains its minimum over all future timelike unit vector $ T_0 $ at
$T_0= \W / \sqrt{ -  \la \W,  \W \ra} $. Now  \eqref{inf-e1}  follows
directly from \eqref{mainest1}.

(iii) At $ \rho = 0 $, we observe that
 $  f $, $   \frac{\p B}{\p \rho} $, $ \frac{\p F}{\p \rho}  $  all
equal $0$  by  \eqref{doff} and \eqref{pB}-\eqref{pF},
hence $ \frac{ \p \tilde{E}}{\p \rho} = 0 $ by \eqref{E(rho)-e2}.
Therefore,  the gradient of $ E(a) $ at $ a = (0, 0, 0)$
is $ - \V$ by \eqref{qle1}. Now if  \eqref{infmly} holds, then
$ a = (0, 0, 0) $ is a critical point of $ E ( a ) $ by
the fact that $ E ( 0 ) = \mly(\Sigma) $, therefore
$ - \V = 0 $. On the other hand, suppose $ - \V = 0$,
by \eqref{mainest1} we have
$ E(\Sigma, X, T_0) \geq \sqrt{ 1 + | a |^2 } \mly ( \Sigma) $,
which implies
$ E(\Sigma, X, T_0) \geq \mly ( \Sigma) $ by
the assumption that $ \mly ( \Sigma) \geq 0 $. Hence,
\eqref{infmly} holds.

This completes the proof of Theorem \ref{estinf}.

\end{proof}

Suppose $ \Sigma $ sits in a time-symmetric slice in $ N$,
then  one would like to compare  the Brown-York mass of $ \Sigma$
 and the Wang-Yau energy of $ \Sigma$. As an
immediate corollary of Theorem \ref{estinf}, we have

\begin{coro} \label{timesym}
Suppose $ \Sigma $ bounds a compact, time-symmetric
hypersurface $ \Omega $ in a spacetime $ N $ satisfying the
dominant energy condition. Suppose $ \Sigma $ has positive
mean curvature $ k $ in $ \Omega $ with respect to the
outward unit normal $ \nu$, then
$$ \inf_{T_0} E( \Sigma, X, T_0 )  =  \mby ( \Sigma ) , $$
and the infimum is achieved at  $ T_0 = ( 1, 0, 0, 0) $.
Moreover, if $ \mby (\Sigma ) > 0 $, then $ (1, 0,0,0)$ is
the unique absolute minimum point of $ E(\Sigma, X, T_0)$
when viewed as a function of $T_0$.
\end{coro}

\begin{proof}
Let $ n $ be the future timelike unit normal to $ \Omega $
in $ N$. Since  $ \Omega $ is time-symmetric  and $ k = -
\la H, \nu \ra > 0$, we have $ H = - k \nu $, $ J = k n $,
and $ \nabla^N_Y \frac{ J}{ | H |}  = \nabla^N_Y n  = 0 $
for any vector $ Y $ tangent to $ \Omega $. Hence,  the
vector field $ V  $ vanishes pointwise on $ \Sigma$. As a
result, $ \V = \int_\Sigma V d v_\Sigma = 0 . $ By Theorem
\ref{estinf} and the fact $ \mly(\Sigma) = \mby(\Sigma)$
under the assumptions, we have \be \label{estuseby}
\begin{split}
  E( \Sigma, X, T_0 )  \ge &  \sqrt{1+ | a |^2}  \mby ( \Sigma) .
   \end{split}
\ee On the other hand, by the positivity results on the
Brown-York mass \cite{ShiTam02}, we have \be \label{pofbym}
 \mby( \Sigma) \geq 0 .
 \ee
Therefore, it follows from \eqref{estuseby} and
\eqref{pofbym} that \be
 E( \Sigma, X, T_0 )  \ge   \mby ( \Sigma) ,
\ee where $ \mby ( \Sigma) $ is the value of $ E(\Sigma, X,
T_0)$ when $ T_0 = ( 1, 0, 0, 0)$. If $ \mby ( \Sigma) > 0
$, \eqref{estuseby} further implies \be
  E( \Sigma, X, T_0 )  >    \mby ( \Sigma)
\ee for any $ T_0  \neq (1, 0, 0, 0) $. The corollary is
thus proved.
\end{proof}

\section{Large sphere limit of  Wang-Yau quasi-local energy  at spatial infinity}

In this section, we study the asymptotical behavior of
$  E( S_r , X_r,T_0) $ in an asymptotically flat spacelike hypersurface
$(M, g, p)$  in a spacetime $ N$, where $ S_r$ is the coordinate sphere
 in a fixed end of $ M$  and $ X_r $  is a suitably chosen embedding of $ S_r $  into $ \R^3 = \{ (0, x) \in \R^{3,1} \}$.

First, we recall the definition of the ADM energy-momentum of $(M, g, p)$. Let $\{ y^i \ | \ i = 1, 2, 3 \}$ be an asymptotic flat coordinate chart on $M$,
 the ADM energy-momentum \cite{adm} of $(M, g, p)$  is a  four covector
$$(E, P_1, P_2,P_3) ,$$
where
\[E=\lim_{r\rightarrow \infty}\frac{1}{16\pi}
\int_{S_r}(\partial_j g_{ij}-\partial_i g_{jj})\nu^i dv_r\]
is the ADM energy of $ (M, g, p)$ and
\[P_k=\lim_{r\rightarrow \infty} \frac{1}{16\pi}\int_{S_r}
2(p_{ik}-\delta_{ik} p_{jj})\nu^i dv_r\] is
 the ADM linear momentum of $(M, g, p)$ in the $ y^k$-direction.
 Here $S_r $ is the coordinate sphere $\{ | y | =r \}$ and
  $ \nu^i \frac{\p }{\p y^i} $ is the outward unit normal to $ S_r $.

Under the assumptions \eqref{af-e1} and \eqref{af-e2}
on $(M, g, p)$, we have the following  fact from \cite[Lemma 2.3]{fst}.
  \begin{lem}\label{fst-l}
  Let $(M,g,p)$ be an asymptotically flat spacelike hypersurface in
  a spacetime $N$. Let $Y=(y^1,y^2,y^3)$ be the asymptotically
  flat coordinates on a fixed end of $ M$.
  There exist an $r_0$ and a constant $C$ independent of $r$ such that for $r\ge r_0$, there is an isometric embedding $ X_r=(x^1,x^2,x^3)$ of $ S_r = \{ | y | =r \}$  into the Euclidean space $\R^3$ such that
  $$
   |X_r-Y|+r ||\nabla X_r-\nabla Y ||_{h_r}+r^2\Big| |H_0|-|H|\Big|+r|n_r^0-n_r|\le C
  $$
  where $X_r$ and $Y$ are considered as $ \R^3$-valued functions on the sphere $S_r$, $H_0$ is the mean curvature vector of $X_r(S_r)$, the gradient $ \nabla (\cdot) $
  and the norm $ || \cdot ||_{h_r} $ are  taken with respect to  the induced
  metric $h_r$ on $S_r$,   $n_r^0$ is the unit outward normal of $X_r(S_r)$, and
 $n_r = Y / | Y | $.
  \end{lem}

We will always work with the embedding $ X_r $ of $ S_r $ provided by Lemma \ref{fst-l}. Given such an $ X_r$, we  let
  \be \label{dfofvrwr}
  \V_r=\V(S_r,X_r), \W_r=\W(S_r,X_r)
  \ee
  which are defined by \eqref{dfofV} and \eqref{dfofW}.

\begin{thm} \label{asympW}
Let $(M,g,p)$, $ S_r $, $Y$ and $X_r$ be as in Lemma
\ref{fst-l}. Let $ \V_r$, $\W_r$ be given in \eqref{dfofvrwr}.
Let $E(S_r,X_r,T_0)$ be the Wang-Yau
quasi-local energy of $ S_r$, where $T_0=(\sqrt{1+|a|^2},a)  $ is any
constant future timelike unit vector in $\R^{3,1}$ with  $a = (a^1, a^2, a^3) $.
Then the followings are true:
\begin{itemize}
  \item [(i)] \be \lim_{r \rightarrow \infty} \W_r = (E, - P_1, -
P_2, - P_3 ), \ee where $(E, -P_1, -P_2, -P_3) $ is
the ADM energy-momentum four vector of $(M, g, p)$.
  \item [(ii)]
  \be\label{asyE-e1}
  E(S_r,X_r,T_0)=-\la T_0, \W_r \ra+\e(r)(1+|a|^2)^\frac12
  \ee
  where $\e(r)$ is a quantity such that $\e(r)\to0$ as
  $r\to\infty$ uniformly in $a$.
  \item[(iii)] Suppose $(E,-P_1,-P_2,-P_3)$ is future timelike.
  Then
  \be \label{inflimitadm}
  \begin{split}
    \lim_{r\to\infty}\inf_{T_0}E(S_r,X_r,T_0) = & \  \madm \\
    =  & \ \lim_{r \to \infty} E(S_r, X_r,T_{0r}^*)   ,
  \end{split}
  \ee
  where $ \madm =\sqrt{E^2-|P|^2}$ is the ADM mass of $(M, g, p)$ and
  $ T_{0r}^* = {\W_r} / {\sqrt{-\la\W_r,\W_r \ra}}$ for sufficiently large $ r$.

\end{itemize}
\end{thm}

\begin{proof}
(i)
By the proof of \cite[Theorem 3.1]{wy3}, we know
 \be \label{wylimit1}
 \lim_{ r \rightarrow \infty}   \frac{1}{8\pi} \int_{S_r} \langle\nabla^{N}_{\nabla\tau} \frac{{J}}{|{H}|},
\frac{{H}}{|{H}|}\rangle d v_r = \sum_{i=1}^3 a^i P_i .
 \ee
By \eqref{cinN}, we have
\begin{equation} \label{cinN2}
\langle\nabla^{N}_{\nabla\tau} \frac{{J}}{|{H}|},
\frac{{H}}{|{H}|}\rangle = - \la a, V \ra ,
\end{equation}
where  $ V=V_r $ is the vector field on $ S_r $ dual to the
$ 1 $-form $ \alpha_{\hat{e}_3} ( \cdot ) $ defined in
\eqref{c1form}. Therefore, by the definition of $ \V_r $ and
\eqref{wylimit1}-\eqref{cinN2},
\be \label{wylimit2}
\lim_{ r \rightarrow \infty}  \la a, \V_r \ra =
  \lim_{ r \rightarrow \infty}   \frac{1}{8\pi} \int_{S_r}
 \la a, V \ra  = - \sum_{i=1}^3 a^i P_i .
\ee
Since $ a = (a^1, a^2, a^3)$ can be chosen arbitrarily,
\eqref{wylimit2} implies
\be
 \lim_{ r \rightarrow \infty}  \V_r = - ( P_1, P_2, P_3).
\ee

Next, we show
\be \label{limofly}
\lim_{ r \rightarrow
\infty} \mly( S_r) = E . \ee By  \cite{fst}, we have $
\lim_{ r \rightarrow \infty}\limits \mby( S_r) = E $. Hence, it
suffices to prove \be \label{dblzero} \lim_{ r \rightarrow
\infty} ( \mly( S_r ) - \mby( S_r) ) = 0. \ee It follows
from the definitions of $ \mly ( S_r ) $ and $ \mby( S_r) $
that
\be \label{difofbl}
\begin{split}
\mly( S_r ) - \mby( S_r) = &
\int_{ S_r} ( k - | H | )   d v_r \\
= & \int_{S_r} ( k - \sqrt{ k^2 - ( \tr_{_{S_r}} p )^2 } )
d v_r ,
\end{split}
\ee where $ k $ is the mean curvature of $ S_r $ in $ (M,
g)$ with respect to the outward unit normal $ \nu_r $ and $
\tr_{_{S_r}} p $ denotes the trace of $ p $ restricted to $
S_r $. Write $ \nu_r = \nu^i \frac{\p}{\p y_i} $, by \cite[Lemma 2.1]{fst}, we have
\be \label{fst-2}
\nu^i = \frac{ y^i}{ r } + O (
r^{-1} ) \ \ \mathrm{and} \ \ k = \frac{2}{r} + O ( r^{-2}).
\ee
Hence, \be
\begin{split}
\tr_{_{S_r}} p =  g^{ij} p_{ij} - p ( \nu_r, \nu_r)
=  O ( r^{-2} ).
\end{split}
\ee Therefore, \be \label{dph}
\begin{split}
 k - \sqrt{ k^2 - ( \tr_{_{S_r}} p )^2 } =
 \frac{( \tr_{_{S_r}} p )^2  }{  k + \sqrt{ k^2 - ( \tr_{_{S_r}} p )^2 }}
 = O ( r^{-3} ) ,
\end{split}
\ee which together with \eqref{difofbl}  implies
\eqref{dblzero}.

(ii) By \eqref{mainest1} in Theorem \ref{estinf}, we have
\be \label{difference}
| E(S_r, X_r, T_0) - ( - \la T_0 , \W_r \ra ) | \leq
C_r  ( 1 + | a|^2)^\frac12,
\ee
where
$$ C_r = \sup_{S_r}
\lf|\lf(\frac{|H_0|^2}{|H|^2}+\frac{|H_0|}{|H|}-2\ri)\ri|
\frac{1}{8\pi} \int_{S_r}  \lf|\lf(|H_0|-|H|\ri)\ri| .
$$
By  Lemma \ref{fst-l} and the fact
$ |H|= 2 r^{-1} + O ( r^{-2}) $ in \eqref{fst-2},
we have
\be \label{limitcr}
 \lim_{r \rightarrow \infty} C_r = 0 .
 \ee
Therefore, \eqref{asyE-e1} follows from \eqref{difference} and \eqref{limitcr}.

(iii)  Since $(E, -P_1, -P_2, -P_3) $ is future timelike,
by (i) $ \W_r $ is future timelike if $ r $ is sufficiently large.
For such an $ r$, by \eqref{inf-e1} in Theorem \ref{estinf},
\begin{equation}\label{inf-e2}
\begin{split}
\sqrt{-\la\W_r,\W_r \ra} + C_r \frac{ \mly(S_r) }{ \sqrt{-\la \W_r,\W_r\ra}  } \ge
& E(S_r, X_r,T_{0r}^*)  \\
\ge &
\inf_{T_0}E(S_r, X_r,T_0)
\ge  \sqrt{-\la \W_r ,\W_r \ra} .
\end{split}
\end{equation}
Now  \eqref{inflimitadm}  follows  from \eqref{inf-e2}, \eqref{limitcr},
\eqref{limofly}
and the fact that
$$ \lim_{r \rightarrow \infty} \sqrt{ - \la \W_r, \W_r \ra } = \madm .$$
This completes the proof of Theorem \ref{asympW}.
\end{proof}

\begin{coro}\label{dominant-t1} Let $(M,g,p)$ be an asymptotically flat spacelike hypersurface in a spacetime
$N$ satisfying the dominant energy condition. With the same notations as in Theorem \ref{asympW}, if
$N$ is not flat along $M$, then
$$
\lim_{r\to\infty}\inf_{ T_0 }E(S_r,X_r,T_0)= \madm.
$$
\end{coro}
\begin{proof} By the positive mass theorem of Schoen-Yau \cite{sy3} and Witten \cite{wi}, the ADM energy-momentum
of $(M, g, p)$ is future timelike unless $N$ is flat along $M$. The result now follows from
Theorem \ref{asympW}.
\end{proof}

Suppose there is an end of $(M, g, p)$ such that $\madm=0$, by \cite{sy3,wi} we know
 $N$ is flat along $M$, $(M,g,p)$ has only one end and can be isometrically embedded in the Minkowski space $ \R^{3,1}$,
 moreover $E=0$, $ P_i=0$, $ i = 1, 2, 3$.
 By \cite{wy2,wy3}, we then have:
$$
\lim_{r\to\infty}\inf_{  \{ \text{admissible $T_0$} \}  }E(S_r,X_r,T_0)=0.
$$
It is still unclear whether the following is true,
$$
\lim_{r\to\infty}\inf_{ T_0 } E(S_r,X_r,T_0)=0.
$$


\end{document}